# STEIN'S METHOD, PALM THEORY AND POISSON PROCESS APPROXIMATION[1]

By Louis H. Y. Chen and Aihua Xia

*National University of Singapore and University of Melbourne*

The framework of Stein's method for Poisson process approximation is presented from the point of view of Palm theory, which is used to construct Stein identities and define local dependence. A general result (Theorem 2.3) in Poisson process approximation is proved by taking the local approach. It is obtained without reference to any particular metric, thereby allowing wider applicability. A Wasserstein pseudometric is introduced for measuring the accuracy of point process approximation. The pseudometric provides a generalization of many metrics used so far, including the total variation distance for random variables and the Wasserstein metric for processes as in Barbour and Brown [*Stochastic Process. Appl.* **43** (1992) 9–31]. Also, through the pseudometric, approximation for certain point processes on a given carrier space is carried out by lifting it to one on a larger space, extending an idea of Arratia, Goldstein and Gordon [*Statist. Sci.* **5** (1990) 403–434]. The error bound in the general result is similar in form to that for Poisson approximation. As it yields the Stein factor $1/\lambda$ as in Poisson approximation, it provides good approximation, particularly in cases where $\lambda$ is large. The general result is applied to a number of problems including Poisson process modeling of rare words in a DNA sequence.

**1. Introduction.** Poisson approximation was developed by Chen (1975) as a discrete version of Stein's normal approximation (1972). It involves the solution of a first-order difference equation, which we call a Stein equation. In extending Poisson approximation to higher dimensions and to Poisson process approximation, Barbour (1988) converted the first-order difference

Received June 2002; revised November 2003.

[1]Supported by Grant R-146-000-013-112 at the National University of Singapore, ARC Discovery Grant DP0209179 from the University of Melbourne and by Grant 01/1/21/19/217 from BMRC of Singapore.

equation into a second-order difference equation and solved it in terms of an immigration-death process. This work was further extended by Barbour and Brown (1992), who introduced a Wasserstein metric on point processes and initiated a program to obtain error bounds of similar order to that on the total variation distance in Poisson approximation. This has been achieved for some special cases by Xia (1997, 2000), and a general result with error bounds of the desired order has been obtained by Brown, Weinberg and Xia (2000).

In this paper, another general result on Poisson process approximation is proved by taking the local approach. It is obtained without reference to any particular metric, thereby allowing wider applicability. In proving this result, the framework of Stein's method is first presented from the point of view of Palm theory, which is used to construct Stein identities and define local dependence. Although the connection between Stein's method and Palm theory has been known to others [see, e.g., Barbour and Månsson (2002)], little of it has been exploited.

In applying the general result, a Wasserstein pseudometric is introduced for measuring the accuracy of point process approximation. The pseudometric provides a generalization of many metrics used so far, including the total variation distance for random variables and the Wasserstein metric for processes as in Barbour and Brown (1992). Also, through the pseudometric, approximation for certain point processes on a given carrier space is carried out by lifting it to one on a larger space, extending an idea of Arratia, Goldstein and Gordon [(1990), Section 3.1], which was refined by Chen [(1998), Section 5].

The error bound in the general result is similar in form to that for Poisson approximation [see, e.g., Arratia, Goldstein and Gordon (1989), Theorem 1]. It is simpler and easier to apply than that in Brown, Weinberg and Xia (2000). As it yields the Stein factor $1/\lambda$ as in Poisson approximation, it provides good approximation, particularly in cases where $\lambda$ is large.

The general result is applied to prove approximation theorems for Matérn hard-core processes and for marked dependent trials. The latter is in turn applied to the classical occupancy problem and rare words in biomolecular sequences. The last application, in fact this paper, is motivated by an interest in modeling the distribution of palindromes in a herpesvirus genome by a Poisson process. In Leung, Choi, Xia and Chen (2002), the Poisson process model is used to provide a mathematical basis for using $r$-scans in determining nonrandom clusters of palindromes in herpesvirus genomes [see also Leung and Yamashita (1999)].

**2. From Palm theory to Stein's method.**    Let $\Gamma$ be a fixed locally compact second countable Hausdorff topological space. Such a space is also a Polish space, that is, a space for which there exists a separable and complete metric in $\Gamma$ which generates the topology. Define $\mathcal{H}$ to be the space of nonnegative integer-valued locally finite measures on $\Gamma$, and let $\mathcal{B}$ be the smallest $\sigma$-algebra in $\mathcal{H}$ making the mappings $\xi \mapsto \xi(C)$ measurable for all relatively compact Borel sets $C \subset \Gamma$. Recall that a point process on $\Gamma$ is a measurable mapping of some fixed probability space into $(\mathcal{H}, \mathcal{B})$ [Kallenberg (1983), page 5]. For a point process $\Xi$ on $\Gamma$ with locally finite mean measure $\boldsymbol{\lambda}$, the point process $\Xi_\alpha$ is said to be a Palm process associated with $\Xi$ at $\alpha \in \Gamma$ if, for any measurable function $f \colon \Gamma \times \mathcal{H} \to \mathbb{R}_+ := [0, \infty)$,

$$(2.1) \qquad \mathbb{E}\left( \int_\Gamma f(\alpha, \Xi) \Xi(d\alpha) \right) = \mathbb{E}\left( \int_\Gamma f(\alpha, \Xi_\alpha) \boldsymbol{\lambda}(d\alpha) \right)$$

[Kallenberg (1983), Chapter 10]. Intuitively,

$$\mathbb{P}(\Xi_\alpha \in B) = \frac{\mathbb{E}[\Xi(d\alpha); \mathbb{1}_{\Xi \in B}]}{\mathbb{E}\Xi(d\alpha)} \qquad \text{for all } B \in \mathcal{B}.$$

An important characterization of Poisson process in the language of Palm theory is that $\Xi$ is a Poisson process if and only if $\mathcal{L}(\Xi_\alpha) = \mathcal{L}(\Xi + \delta_\alpha)$ $\boldsymbol{\lambda}$-a.s., where $\delta_\alpha$ is the Dirac measure at $\alpha$. This highlights an idea of Poisson process approximation: if we define

$$Df(\xi) := \int_\Gamma f(x, \xi)\xi(dx) - \int_\Gamma f(x, \xi + \delta_x)\boldsymbol{\lambda}(dx),$$

then $\mathcal{L}(\Xi)$ is close to the Poisson process distribution over $\Gamma$ with mean measure $\boldsymbol{\lambda}$, denoted as $\mathrm{Po}(\boldsymbol{\lambda})$, in terms of a certain metric if, for the set of suitable corresponding test functions $f \colon \Gamma \times \mathcal{H} \to \mathbb{R} := (-\infty, \infty)$,

$$(2.2) \qquad \mathbb{E}Df(\Xi) \approx 0.$$

In other words, for a function $g \colon \mathcal{H} \to \mathbb{R}$, if we can find a solution $f_g$ to the equation

$$(2.3) \qquad g(\xi) - \mathrm{Po}(\boldsymbol{\lambda})(g) = Df(\xi),$$

then the distance between the distribution of $\Xi$ and $\mathrm{Po}(\boldsymbol{\lambda})$ is achieved by the supremum of $|\mathbb{E}Df_g(\Xi)|$ over the class of $g$ which defines the metric. Equation (2.3) is known as a Stein equation. If there exists a function $h \colon \mathcal{H} \to \mathbb{R}$ such that $f(x, \xi) = h(\xi - \delta_x) - h(\xi)$, then

$$Df(\xi) = \int_\Gamma [h(\xi + \delta_x) - h(\xi)]\boldsymbol{\lambda}(dx) + \int_\Gamma [h(\xi - \delta_x) - h(\xi)]\xi(dx) := \mathcal{A}h(\xi).$$

It is known that $\mathcal{A}$ is the generator of an $\mathcal{H}$-valued immigration-death process $\mathbf{Z}_\xi(t)$ with immigration intensity $\boldsymbol{\lambda}$ and unit per capita death rate, where $\mathbf{Z}_\xi(0) = \xi$. This fact was noted by Barbour (1988), who developed a probabilistic approach to Stein's method for multivariate Poisson and Poisson process approximations. The equilibrium distribution of $\mathbf{Z}_\xi$ is a Poisson process with mean measure $\boldsymbol{\lambda}$. The idea of introducing a Markov point process is to exploit the probabilistic properties of the Markov process for obtaining bounds on the metrics of interest [see Barbour and Brown (1992) and Brown and Xia (2000)].

For $\xi \in \mathcal{H}$ and a Borel set $B \subset \Gamma$, we define $\xi|_B$ as the restriction of $\xi$ to $B$, that is, $\xi|_B(C) = \xi(B \cap C)$ for Borel sets $C \subset \Gamma$. Let $\Xi$ be a point process on $\Gamma$ with Palm processes $\{\Xi_\alpha\}$. Assume that for each $\alpha$ there is a Borel set $A_\alpha \subset \Gamma$ such that $\alpha \in A_\alpha$ and the mapping

$$(2.4) \qquad \Gamma \times \mathcal{H} \to \Gamma \times \mathcal{H} : (\alpha, \xi) \mapsto (\alpha, \xi^{(\alpha)})$$

is product measurable, where $\xi^{(\alpha)} := \xi|_{A_\alpha^c}$. Note that $\xi^{(\alpha)}$ does not refer to the Palm measure. As the measurability of (2.4) is often hard to check, we give a sufficient condition for (2.4) to hold: $A = \{(x, y) : y \in A_x, x \in \Gamma\}$ is a measurable set of the product space $\Gamma^2 := \Gamma \times \Gamma$. We give a brief proof for the sufficiency. By the monotone class theorem, it suffices to show that the mapping $M_A(\alpha, \xi) := (\alpha, \xi^{(\alpha)})$ is measurable for rectangular sets $A = B_1 \times B_2$, where $B_1$ and $B_2$ are measurable subsets of $\Gamma$. Indeed,

$$M_{B_1 \times B_2}(x, \xi) = \begin{cases} (x, \xi|_{B_2^c}), & \text{if } x \in B_1, \\ (x, \xi), & \text{if } x \notin B_1, \end{cases}$$

is measurable.

The requirement of $A$ being measurable in $\Gamma^2$ is almost necessary. To see this, let $\Gamma = [0,1]$, $A = B_1 \times B_2$, where $B_1 \subset \Gamma$ is not Borel measurable [Nielsen (1997), page 128, 9.16(h)] and $B_2 \subset \Gamma$ is a Borel set. Define $C = \{\xi : \xi(B_2) \neq 0\}$; then $M_A^{-1}(\Gamma \times C) = B_1^c \times C$ is not a measurable set of $\Gamma \times \mathcal{H}$.

REMARK 2.1.    In Barbour and Brown [(1992), page 15], it is proved that if $A_\alpha$ is a ball of fixed radius, then the mapping in (2.4) is measurable.

We define $\Xi$ to be *locally dependent with neighborhoods* $(A_\alpha; \alpha \in \Gamma)$ if

$$\mathcal{L}((\Xi_\alpha)^{(\alpha)}) = \mathcal{L}(\Xi^{(\alpha)}), \qquad \lambda\text{-a.s.}$$

LEMMA 2.2.    *The following statements are equivalent:*

(a) $\mathbb{E} \int_\Gamma f(\alpha, \Xi^{(\alpha)} + \delta_\alpha) \Xi(d\alpha) = \mathbb{E} \int_\Gamma f(\alpha, \Xi^{(\alpha)} + \delta_\alpha) \lambda(d\alpha)$ *for all measurable* $f : \Gamma \times \mathcal{H} \to \mathbb{R}_+$.
(b) $\mathcal{L}((\Xi_\alpha)^{(\alpha)}) = \mathcal{L}(\Xi^{(\alpha)})$, $\lambda$-*a.s.*

PROOF.    By the definition of Palm process, we have

$$(2.5) \qquad \mathbb{E} \int_\Gamma f(\alpha, \Xi^{(\alpha)} + \delta_\alpha) \Xi(d\alpha) = \mathbb{E} \int_\Gamma f(\alpha, (\Xi_\alpha)^{(\alpha)} + \delta_\alpha) \lambda(d\alpha).$$

Hence, (b) implies (a). Now assume (a). With the vague topology, $\mathcal{H}$ is a Polish space [see Kallenberg (1983), page 95], so there exists a sequence of bounded uniformly continuous functions $(f_j; j \geq 1)$ on $\mathcal{H}$ which form a determining class [Billingsley (1968), page 15]: for every two probability measures $\mathbf{Q}_1$ and $\mathbf{Q}_2$ on $\mathcal{H}$, if $\int f_j \, d\mathbf{Q}_1 = \int f_j \, d\mathbf{Q}_2$ for all $j \geq 1$, then $\mathbf{Q}_1 = \mathbf{Q}_2$ [see Parthasarathy (1967), Theorem 6.6]. By taking $f(\alpha, \xi + \delta_\alpha) = k(\alpha)f_j(\xi)$, it follows from (2.5) that

$$\int_\Gamma k(\alpha)[\mathbb{E}f_j(\Xi^{(\alpha)})]\boldsymbol{\lambda}(d\alpha) = \int_\Gamma k(\alpha)[\mathbb{E}f_j((\Xi_\alpha)^{(\alpha)})]\boldsymbol{\lambda}(d\alpha)$$

for all bounded measurable functions $k \colon \Gamma \to \mathbb{R}_+$ and $f_j$. Fixing $f_j$ and allowing $k$ to vary, we have $\mathbb{E}f_j(\Xi^{(\alpha)}) = \mathbb{E}f_j((\Xi_\alpha)^{(\alpha)})$, $\boldsymbol{\lambda}$-a.s. Now vary $f_j$ and (b) follows. □

In general, a point process is not necessarily locally dependent, but Lemma 2.2 suggests that, in a loose sense,

(2.6) $$\mathcal{L}((\Xi_\alpha)^{(\alpha)}) \approx \mathcal{L}(\Xi^{(\alpha)}), \qquad \boldsymbol{\lambda}\text{-a.s.}$$

if and only if

(2.7) $$\begin{aligned} &\mathbb{E}\int_\Gamma f(\alpha, \Xi^{(\alpha)} + \delta_\alpha)\Xi(d\alpha) \\ &\qquad \approx \mathbb{E}\int_\Gamma f(\alpha, \Xi^{(\alpha)} + \delta_\alpha)\boldsymbol{\lambda}(d\alpha) \qquad \text{for suitable } f. \end{aligned}$$

This will be our guiding principle in proving Theorem 2.3 using the local approach, as follows [an extension of the approach of Chen (1975) which was elaborated by Barbour and Brown (1992)]:

$$\begin{aligned} &\mathbb{E}\int_\Gamma f(\alpha, \Xi)\Xi(d\alpha) \\ &\qquad = \mathbb{E}\int_\Gamma [f(\alpha, \Xi) - f(\alpha, \Xi^{(\alpha)} + \delta_\alpha)]\Xi(d\alpha) \\ &\qquad\quad + \mathbb{E}\int_\Gamma f(\alpha, \Xi^{(\alpha)} + \delta_\alpha)[\Xi(d\alpha) - \boldsymbol{\lambda}(d\alpha)] \\ &\qquad\quad + \mathbb{E}\int_\Gamma [f(\alpha, \Xi^{(\alpha)} + \delta_\alpha) - f(\alpha, \Xi + \delta_\alpha)]\boldsymbol{\lambda}(d\alpha) \\ &\qquad\quad + \mathbb{E}\int_\Gamma f(\alpha, \Xi + \delta_\alpha)\boldsymbol{\lambda}(d\alpha), \end{aligned}$$

which implies

$$\mathbb{E}Df(\Xi) = \mathbb{E}\int_\Gamma [f(\alpha, \Xi) - f(\alpha, \Xi^{(\alpha)} + \delta_\alpha)]\Xi(d\alpha)$$

$$(2.8) \qquad + \mathbb{E} \int_\Gamma f(\alpha, \Xi^{(\alpha)} + \delta_\alpha)[\Xi(d\alpha) - \boldsymbol{\lambda}(d\alpha)]$$

$$+ \mathbb{E} \int_\Gamma [f(\alpha, \Xi^{(\alpha)} + \delta_\alpha) - f(\alpha, \Xi + \delta_\alpha)]\boldsymbol{\lambda}(d\alpha).$$

Hence, a bound on $\mathbb{E} D f_g(\Xi)$ can be obtained by bounding the right-hand side of (2.8).

There are two ways to handle the second term in (2.8): one uses coupling and the other involves Janossy densities [Janossy (1950) and Daley and Vere-Jones (1988)]. For a finite point process $\Xi$, that is $\mathbb{P}(|\Xi| < \infty) = 1$, there exist measures $(J_n)_{n \geq 1}$ such that, for measurable functions $f : \mathcal{H} \to \mathbb{R}_+$,

$$\mathbb{E} f(\Xi) = \sum_{n \geq 0} \int_{\Gamma^n} \frac{1}{n!} f\left(\sum_{i=1}^n \delta_{x_i}\right) J_n(dx_1, \ldots, dx_n).$$

The term $(n!)^{-1} J_n(dx_1, \ldots, dx_n)$ can be intuitively explained as the probability of $\Xi$ having $n$ points and these points being located near $(x_1, \ldots, x_n)$. The measures $(J_n)_{n \geq 1}$ are called *Janossy measures* by Srinivasan (1969).

Suppose there is a reference measure $\boldsymbol{\nu}$ on $\Gamma$ such that, for each $n \geq 1$, $J_n$ is absolutely continuous with respect to $\boldsymbol{\nu}^n$. Then, by the Radon–Nikodym theorem, the derivatives $j_n$ of $J_n$ with respect to $\boldsymbol{\nu}^n$ exist, so that

$$\mathbb{E} f(\Xi) = \sum_{n \geq 0} \int_{\Gamma^n} \frac{1}{n!} f\left(\sum_{i=1}^n \delta_{x_i}\right) j_n(x_1, \ldots, x_n)\boldsymbol{\nu}^n(dx_1, \ldots, dx_n).$$

The derivatives $(j_n)_{n \geq 1}$ are called *Janossy densities*.

The density of the mean measure $\boldsymbol{\lambda}$ of a finite point process $\Xi$ with respect to $\boldsymbol{\nu}$ can be expressed by its Janossy densities $(j_n)_{n \geq 1}$ as

$$\phi(x) = \sum_{m \geq 0} \int_{\Gamma^m} \frac{1}{m!} j_{m+1}(x, x_1, \ldots, x_m)\boldsymbol{\nu}^m(dx_1, \ldots, dx_m),$$

where the term with $m = 0$ is interpreted as $j_1(x)$ [Daley and Vere-Jones (1988), page 133].

When the point process is simple, the Janossy densities can also be used to describe the conditional probability density of a point being at $\alpha$, given the configuration $\Xi^{(\alpha)}$ of $\Xi$ outside $A_\alpha$. More precisely, let $m \in \mathbf{N}$ be fixed and $\boldsymbol{\beta} = (\beta_1, \ldots, \beta_m) \in (A_\alpha^c)^m$, and define

$$(2.9) \qquad \mathcal{G}(\alpha, \boldsymbol{\beta}) := \frac{\sum_{r \geq 0} \int_{A_\alpha^r} j_{m+r+1}(\alpha, \boldsymbol{\beta}, \boldsymbol{\gamma})(r!)^{-1} \boldsymbol{\nu}^r(d\boldsymbol{\gamma})}{\sum_{s \geq 0} \int_{A_\alpha^s} j_{m+s}(\boldsymbol{\beta}, \boldsymbol{\eta})(s!)^{-1} \boldsymbol{\nu}^s(d\boldsymbol{\eta})},$$

where the term with $r = 0$ is interpreted as $j_{m+1}(\alpha, \boldsymbol{\beta})$ and the term with $s = 0$ as $j_m(\boldsymbol{\beta})$. Then $\mathcal{G}(\alpha, \boldsymbol{\beta})$ is the conditional density of a point being near $\alpha$ given that $\Xi^{(\alpha)}$ is $\sum_{i=1}^{m} \delta_{\beta_i}$. Direct verification gives that, for any bounded measurable function $f$ over $\mathcal{H}$,

$$(2.10) \quad \mathbb{E}\left(\int_\Gamma f(\alpha, \Xi^{(\alpha)}) \Xi(d\alpha)\right) = \mathbb{E}\left(\int_\Gamma f(\alpha, \Xi^{(\alpha)}) \mathcal{G}(\alpha, \Xi^{(\alpha)}) \boldsymbol{\nu}(d\alpha)\right).$$

For each $f : \Gamma \times \mathcal{H} \to \mathbb{R}_+$, $\xi \in \mathcal{H}$, write $\xi(A_x) = m$ and define

$$(\delta f)(x, \xi) = \sup_{\{z_1, \ldots, z_m\} \subset \Gamma} \max_{0 \le j \le m-1} \left| f\left(x, \xi^{(x)} + \delta_x + \sum_{i=1}^{j} \delta_{z_i}\right) \right.$$
$$\left. - f\left(x, \xi^{(x)} + \delta_x + \sum_{i=1}^{j+1} \delta_{z_i}\right) \right|,$$

where the right-hand side is interpreted as 0 if $m = 0$. Combining (2.3) and (2.8) gives:

THEOREM 2.3. *For each bounded measurable function* $g : \mathcal{H} \to \mathbb{R}_+$,

$$|\mathbb{E}g(\Xi) - \mathrm{Po}(\boldsymbol{\lambda})(g)|$$
$$(2.11) \qquad \le \mathbb{E}\int_{\alpha \in \Gamma}(\delta f_g)(\alpha, \Xi)(\Xi(A_\alpha) - 1)\Xi(d\alpha) + \min\{\varepsilon_1(g, \Xi), \varepsilon_2(g, \Xi)\}$$
$$\qquad + \mathbb{E}\int_{\alpha \in \Gamma}(\delta f_g)(\alpha, \Xi)\boldsymbol{\lambda}(d\alpha)\Xi(A_\alpha),$$

*where*

$$(2.12) \quad \varepsilon_1(g, \Xi) = \mathbb{E}\int_{\alpha \in \Gamma}|f_g(\alpha, \Xi^{(\alpha)} + \delta_\alpha)||\mathcal{G}(\alpha, \Xi^{(\alpha)}) - \phi(\alpha)|\boldsymbol{\nu}(d\alpha),$$

*which is valid if* $\Xi$ *is a simple point process, and*

$$(2.13) \quad \varepsilon_2(g, \Xi) = \mathbb{E}\int_{\alpha \in \Gamma}|f_g(\alpha, \Xi^{(\alpha)} + \delta_\alpha) - f_g(\alpha, (\Xi_\alpha)^{(\alpha)} + \delta_\alpha)|\boldsymbol{\lambda}(d\alpha).$$

REMARK 2.4. How judicious $(A_\alpha; \alpha \in \Gamma)$ are chosen is reflected in the upper bound in (2.11), and (2.13) suggests that $(A_\alpha; \alpha \in \Gamma)$ should normally be chosen such that (2.6) holds.

**3. Poisson process approximation in Wasserstein pseudometric.** We now look at special test functions $g$ which define metrics of our interest. We begin with a pseudometric $\rho_0$ on $\Gamma$ bounded by 1 [cf. Barbour and Brown (1992)]. In order for Theorem 2.3 to be applicable, we assume that the topology generated by $\rho_0$ is weaker than the given topology of $\Gamma$. Let $\mathcal{K}$ stand for the set of $\rho_0$-Lipschitz functions $k\colon \Gamma \to [-1,1]$ such that $|k(\alpha) - k(\beta)| \le \rho_0(\alpha,\beta)$ for all $\alpha, \beta \in \Gamma$. The first Wasserstein pseudometric $\rho_1$ is defined on $\mathcal{H}$ by

$$\rho_1(\xi_1,\xi_2) = \begin{cases} 1, & \text{if } |\xi_1| \ne |\xi_2|, \\ |\xi_1|^{-1} \sup_{k\in\mathcal{K}} \left| \int k\,d\xi_1 - \int k\,d\xi_2 \right|, & \text{if } |\xi_1| = |\xi_2| > 0, \end{cases}$$

where $|\xi_i|$ is the total mass of $\xi_i$. A pseudometric $\rho_1''$ equivalent to $\rho_1$ can be defined as follows [cf. Brown and Xia (1995)]: for two configurations $\xi_1 = \sum_{i=1}^n \delta_{y_i}$ and $\xi_2 = \sum_{i=1}^m \delta_{z_i}$ with $m \ge n$,

$$\rho_1''(\xi_1,\xi_2) = \min_{\pi} \sum_{i=1}^n \rho_0(y_i, z_{\pi(i)}) + (m - n),$$

where $\pi$ ranges all permutations of $(1,\dots,m)$.

Let $\mathcal{F}$ denote the set of $\rho_1$-Lipschitz functions on $\mathcal{H}$ such that $|f(\xi_1) - f(\xi_2)| \le \rho_1(\xi_1,\xi_2)$ for all $\xi_1$ and $\xi_2 \in \mathcal{H}$. The second Wasserstein pseudometric is defined on probability measures on $\mathcal{H}$ with respect to $\rho_1$ by

$$\rho_2(\mathbf{Q}_1,\mathbf{Q}_2) = \sup_{f\in\mathcal{F}} \left| \int f\,d\mathbf{Q}_1 - \int f\,d\mathbf{Q}_2 \right|.$$

The use of a pseudometric $\rho_0$ provides not only generality but also wider applicability. For example, if we choose $\rho_0(x,y) \equiv 0$, then

$$\rho_2(\mathbf{Q}_1,\mathbf{Q}_2) = d_{\mathrm{TV}}(\mathcal{L}(|X_1|), \mathcal{L}(|X_2|)),$$

the total variation distance between $\mathcal{L}(|X_1|)$ and $\mathcal{L}(|X_2|)$, where $X_i$ has distribution $\mathbf{Q}_i$, $i = 1,2$. It is known that, for $g = \mathbb{1}_B$ with $B \subset \mathbf{Z}_+ := \{0,1,2,\dots\}$,

$$(\delta f_g)(x,\xi) \le \frac{1 - e^{-\lambda}}{\lambda}, \qquad |f_g| \le 1 \wedge \sqrt{\frac{2}{e\lambda}},$$

where, and throughout this paper, $\lambda$ is the total mass of $\boldsymbol{\lambda}$ and is assumed to be finite [see Barbour, Holst and Janson (1992) and Brown and Xia (2001)]. So Theorem 2.3 gives:

THEOREM 3.1. *We have*

$$d_{\mathrm{TV}}(\mathcal{L}(\Xi(\Gamma)), \mathrm{Po}(\lambda)) \le \frac{1 - e^{-\lambda}}{\lambda} \mathbb{E} \int_{\alpha\in\Gamma} (\Xi(A_\alpha) - 1)\Xi(d\alpha) + \min\{\varepsilon_1, \varepsilon_2\}$$
$$+ \frac{1 - e^{-\lambda}}{\lambda} \int_{\alpha\in\Gamma} \boldsymbol{\lambda}(A_\alpha)\boldsymbol{\lambda}(d\alpha),$$

*where*

$$\varepsilon_1 = 1 \wedge \sqrt{\frac{2}{e\lambda}} \int_{\alpha\in\Gamma} \mathbb{E}|\mathcal{G}(\alpha, \Xi^{(\alpha)}) - \phi(\alpha)|\boldsymbol{\nu}(d\alpha),$$

*which is valid for $\Xi$ simple, and*

$$\varepsilon_2 = \frac{1-e^{-\lambda}}{\lambda}\int_{\alpha\in\Gamma}\mathbb{E}|\,|\Xi^{(\alpha)}| - |(\Xi_\alpha)^{(\alpha)}|\,|\,\boldsymbol{\lambda}(d\alpha).$$

Theorem 3.1 with $\varepsilon_1$ is a generalization of Chen (1975) [see also Barbour and Brown (1992)] and with $\varepsilon_2$ allows the use of the coupling approach [see Barbour and Brown (1992)].

Another example is in Section 4, where it is possible to introduce an index space so that the results also include the approximation in distribution by a Poisson process to discrete sums of the form $\sum_{i=1}^{n} X_i \delta_{Y_i}$, where $Y_i$ is a random mark associated with $X_i$, as in Arratia, Goldstein and Gordon (1989).

We now establish a general statement of this section. As the arguments in Barbour and Brown (1992) and Brown and Xia (2001) never rely on the property that $\rho_0(x,y)=0$ implies $x=y$, the results are still valid for $\rho_0$ and the pseudometrics $\rho_1$ and $\rho_2$ generated from $\rho_0$. The following two lemmas are taken from Barbour and Brown (1992) and Brown and Xia (2001).

LEMMA 3.2. *For each $\rho_1$-Lipschitz function $g \in \mathcal{F}$, $x,y \in \Gamma$ and $\xi \in \mathcal{H}$ with $|\xi| = n$, the solution $f_g$ of (2.3) satisfies*

$$(3.1) \qquad |f_g(x,\xi+\delta_x+\delta_y) - f_g(x,\xi+\delta_x)| \leq \frac{5}{\lambda} + \frac{3}{n+1},$$

$$(3.2) \qquad |f_g(y,\xi+\delta_y)| \leq 1 \wedge 1.65\lambda^{-1/2}.$$

LEMMA 3.3. *For each $g \in \mathcal{F}$, $\xi,\eta \in \mathcal{H}$ and $x \in \Gamma$,*

$$|f_g(x,\xi+\delta_x) - f_g(x,\eta+\delta_x)|$$
$$\leq \frac{2}{|\eta| \wedge |\xi| + 1}[\rho_1''(\xi,\eta) - ||\eta| - |\xi||] + \left(\frac{5}{\lambda} + \frac{3}{|\eta| \wedge |\xi| + 1}\right)||\eta| - |\xi||$$
$$\leq \left(\frac{5}{\lambda} + \frac{3}{|\eta| \wedge |\xi| + 1}\right)\rho_1''(\xi,\eta).$$

With the above two lemmas, we write another version of Theorem 2.3.

THEOREM 3.4. *We have*

$$\rho_2(\mathcal{L}\Xi, \mathrm{Po}(\boldsymbol{\lambda}))$$
$$(3.3) \qquad \leq \mathbb{E}\int_{\alpha\in\Gamma}\left(\frac{5}{\lambda} + \frac{3}{|\Xi^{(\alpha)}|+1}\right)(\Xi(A_\alpha)-1)\Xi(d\alpha) + \min\{\varepsilon_1,\varepsilon_2\}$$
$$\qquad + \mathbb{E}\int_{\alpha\in\Gamma}\int_{\beta\in A_\alpha}\left(\frac{5}{\lambda} + \frac{3}{|(\Xi_\beta)^{(\alpha)}|+1}\right)\boldsymbol{\lambda}(d\alpha)\boldsymbol{\lambda}(d\beta),$$

*where*

$$(3.4) \quad \varepsilon_1 = \left(1 \wedge (1.65\lambda^{-1/2})\right) \int_{\alpha \in \Gamma} \mathbb{E}|\mathcal{G}(\alpha, \Xi^{(\alpha)}) - \phi(\alpha)|\boldsymbol{\nu}(d\alpha),$$

$$(3.5) \quad \varepsilon_2 = \mathbb{E} \int_{\alpha \in \Gamma} \left(\frac{5}{\lambda} + \frac{3}{|(\Xi_\alpha)^{(\alpha)}| \wedge |\Xi^{(\alpha)}| + 1}\right) \rho_1''((\Xi_\alpha)^{(\alpha)}, \Xi^{(\alpha)})\boldsymbol{\lambda}(d\alpha).$$

In many applications, we can obtain the Stein factor $1/\lambda$ from the terms

$$(|\Xi^{(\alpha)}| + 1)^{-1}, \qquad (|(\Xi_\alpha)^{(\alpha)}| \wedge |\Xi^{(\alpha)}| + 1)^{-1}, \qquad (|(\Xi_\beta)^{(\alpha)}| + 1)^{-1},$$

by applying Lemma 3.5.

LEMMA 3.5 [Brown, Weinberg and Xia ([2000](#)), Lemma 3.1]. *For a random variable $X \geq 1$,*

$$\mathbb{E}\left(\frac{1}{X}\right) \leq \frac{\sqrt{\kappa(1 + \kappa/4)} + 1 + \kappa/2}{\mathbb{E}(X)},$$

*where $\kappa = \text{Var}(X)/\mathbb{E}(X)$.*

COROLLARY 3.6. *If $\Xi$ is a locally dependent point process with neighborhoods $(A_\alpha; \alpha \in \Gamma)$, then*

$$(3.6)
\begin{aligned}
\rho_2(\mathcal{L}\Xi, \text{Po}(\boldsymbol{\lambda})) &\leq \mathbb{E} \int_{\alpha \in \Gamma} \left(\frac{5}{\lambda} + \frac{3}{|\Xi^{(\alpha)}| + 1}\right)(\Xi(A_\alpha) - 1)\Xi(d\alpha) \\
&\quad + \int_{\alpha \in \Gamma} \int_{\beta \in A_\alpha} \left(\frac{5}{\lambda} + \mathbb{E}\frac{3}{|\Xi^{(\alpha\beta)}| + 1}\right)\boldsymbol{\lambda}(d\alpha)\boldsymbol{\lambda}(d\beta),
\end{aligned}$$

*where $\xi^{(\alpha\beta)} = \xi|_{A_\alpha^c \cap A_\beta^c}$.*

REMARK 3.7. Since

$$\int_{\alpha \in \Gamma} \frac{\Xi(A_\alpha)}{|\Xi^{(\alpha)}| + 1}\Xi(d\alpha) = \int_{\alpha \in \Gamma} \int_{\beta \in A_\alpha} \frac{1}{|\Xi^{(\alpha)}| + 1}\Xi(d\beta)\Xi(d\alpha)$$

$$\leq \int_{\alpha \in \Gamma} \int_{\beta \in A_\alpha} \frac{1}{|\Xi^{(\alpha\beta)}| + 1}\Xi(d\beta)\Xi(d\alpha),$$

to simplify the first term of (3.6) using the assumption of local dependence, it is tempting to ask whether

$$\mathbb{E}\frac{1}{|\Xi^{(\alpha\beta)}| + 1}\Xi(d\beta)\Xi(d\alpha) = \mathbb{E}\frac{1}{|\Xi^{(\alpha\beta)}| + 1}\mathbb{E}\Xi(d\beta)\Xi(d\alpha).$$

The answer is generally negative, although it might be true in many applications, as shown in Section 5. To see this, let $\mathbb{P}(B_i) = q = 0.1$ for $i = 1, 2, 3$, $\mathbb{P}(B_iB_j) = q^2$ for $1 \leq i \neq j \leq 3$ and $\mathbb{P}(B_1B_2B_3) = 2q^3$. Set $\Gamma = \{1, 2, 3\}$, $\Xi(\{i\}) = \mathbb{1}_{B_i}$, $1 \le i \le 3$, and $A_1 = A_2 = \{1, 2\}$ and $A_3 = \{1, 3\}$; then $\Xi$ is locally dependent with neighborhoods $(A_i; i \in \Gamma)$. However, direct calculation gives

$$\mathbb{E}\frac{1}{\Xi(\{3\}) + 1}\Xi(\{1\})\Xi(\{2\}) = q^2 - q^3$$

and

$$\mathbb{E}\frac{1}{\Xi(\{3\}) + 1}\mathbb{E}\Xi(\{1\})\Xi(\{2\}) = (1 - 0.5q)q^2,$$

so

$$\mathbb{E}\frac{1}{\Xi(\{3\}) + 1}\Xi(\{1\})\Xi(\{2\}) \ne \mathbb{E}\frac{1}{\Xi(\{3\}) + 1}\mathbb{E}\Xi(\{1\})\Xi(\{2\}).$$

**4. Sums of marked dependent trials.** The case of Poisson process approximation for sums of marked dependent trials is of particular interest as it has applications in computational biology, occupancy and random graphs. We devote this section to this case.

Let $I_i$, $i \in \mathcal{I}$, be dependent indicators with $\mathcal{I}$ a finite or infinitely countable index space and

$$\mathbb{P}(I_i = 1) = 1 - \mathbb{P}(I_i = 0) = p_i, \qquad i \in \mathcal{I}.$$

Let $\mathcal{U}_i$, $i \in \mathcal{I}$, be $\mathcal{S}$-valued independent random elements, where $\mathcal{S}$ is a locally compact second countable Hausdorff space with metric $d_0$ bounded by 1. Assume that $\{\mathcal{U}_i, \ i \in \mathcal{I}\}$ is independent of $\{I_i, \ i \in \mathcal{I}\}$. Our interest is to approximate the distribution of $\mathcal{M} := \sum_{i \in \mathcal{I}} I_i \delta_{\mathcal{U}_i}$ by that of a Poisson process.

Let $\mathcal{H}(\mathcal{S})$ be the space of nonnegative integer-valued locally finite measures on $\mathcal{S}$. The metric $d_0$ will generate the first Wasserstein metric $d_1$ on $\mathcal{H}(\mathcal{S})$ and second Wasserstein metric $d_2$ on probability measures on $\mathcal{H}(\mathcal{S})$ as in Section 3 [see also Barbour and Brown (1992)]. For each $i \in \mathcal{I}$, let $A_i \subset \mathcal{I}$ such that $i \in A_i$. Let $\mu_i = \mathcal{L}(\mathcal{U}_i)$, the law of $\mathcal{U}_i$, $i \in \mathcal{I}$, and let $\boldsymbol{\lambda} = \sum_{i \in \mathcal{I}} p_i \mu_i$. Define $V_i = \sum_{j \notin A_i} I_j$.

THEOREM 4.1. *We have* $\lambda = \sum_{i \in \mathcal{I}} p_i$ *and*

$$
\begin{aligned}
d_2(\mathcal{L}\mathcal{M}, \mathrm{Po}(\boldsymbol{\lambda})) \le {} & \mathbb{E}\sum_{i \in \mathcal{I}} \sum_{j \in A_i \setminus \{i\}} \left(\frac{5}{\lambda} + \frac{3}{V_i + 1}\right) I_i I_j + \min\{\varepsilon_1, \varepsilon_2\} \\
& + \sum_{i \in \mathcal{I}} \sum_{j \in A_i} \left(\frac{5}{\lambda} + \mathbb{E}\left[\frac{3}{V_i + 1}\Big| I_j = 1\right]\right) p_i p_j,
\end{aligned}
$$
(4.1)

*where*

$$\varepsilon_1 = (1 \wedge 1.65\lambda^{-1/2}) \sum_{i \in \mathcal{I}} \mathbb{E}|\mathbb{E}(I_i|I_j, j \notin A_i) - p_i|,$$

$$\varepsilon_2 = \mathbb{E}\sum_{i \in \mathcal{I}}\left(\frac{5}{\lambda} + \frac{3}{V_i \wedge \sum_{j \notin A_i} J_{ji} + 1}\right) \sum_{j \notin A_i} |J_{ji} - I_j|p_i,$$

*and* $(J_{ji}; j \in \mathcal{I})$ *and* $(I_j; j \in \mathcal{I})$ *are defined on the same probability space with*

$$\mathcal{L}(J_{ji}; j \in \mathcal{I}) = \mathcal{L}(I_j; j \in \mathcal{I}|I_i = 1).$$

REMARK 4.2. The bound in (4.1) does not depend on the distribution of the marks $(\mathcal{U}_i)_{i \in \mathcal{I}}$, since the mean measure of the approximating Poisson process has been chosen to reflect the contribution of the marks.

REMARK 4.3. Since $\mathcal{M}$ is in general not a *simple* point process, the Janossy density approach via (2.9) is not applicable. Also, due to the structure of $\mathcal{M}$, the neighborhoods $\{A_\alpha, \alpha \in \mathcal{S}\}$ cannot be determined. By introducing a pseudometric and by lifting the process $\mathcal{M}$ from $\mathcal{S}$ to a larger carrier space $\Gamma = \mathcal{S} \times \mathcal{I}$, the lifted process becomes *simple* and the neighborhoods $\{A_\alpha, \alpha \in \Gamma\}$ determinable.

PROOF OF THEOREM 4.1. We consider the approximation on the lifted space $\Gamma = \mathcal{S} \times \mathcal{I}$ with pseudometric $\rho_0((s,i),(t,j)) = d_0(s,t)$. For each $\xi_l \in \mathcal{H}(\Gamma)$ ($l$ means lifted), define $\xi \in \mathcal{H}(\mathcal{S})$ by $\xi(ds) = \sum_{i \in \mathcal{I}} \xi_l(ds, \{i\})$. Let $\mathcal{M}_l(ds, \{i\}) = I_i\delta_{\mathcal{U}_i}(ds)$ and let $\boldsymbol{\lambda}_l(ds, \{i\}) = p_i\mu_i(ds)$. Then $\mathcal{M}_l$ is a simple point process on $\Gamma$, $\mathcal{M}(ds) = \sum_{i \in \mathcal{I}} \mathcal{M}_l(ds, \{i\})$, $\boldsymbol{\lambda}(ds) = \sum_{i \in \mathcal{I}} \boldsymbol{\lambda}_l(ds, \{i\})$, and

$$\rho_2(\mathcal{L}\mathcal{M}_l, \mathrm{Po}(\boldsymbol{\lambda}_l)) = d_2(\mathcal{L}\mathcal{M}, \mathrm{Po}(\boldsymbol{\lambda})).$$

For each $(s,i) \in \Gamma$, define $A_{(s,i)} := \mathcal{S} \times A_i$. Then $|\mathcal{M}_l^{((s,i))}| = V_i$.

The first term in the upper bound of (3.3) becomes

$$\mathbb{E}\int_{(s,i) \in \Gamma}\left(\frac{5}{\lambda} + \frac{3}{V_i + 1}\right)(\mathcal{M}_l(A_{(s,i)}) - 1)I_i\delta_{\mathcal{U}_i}(ds)$$

$$= \mathbb{E}\sum_{i \in \mathcal{I}}\left(\frac{5}{\lambda} + \frac{3}{V_i + 1}\right)\left(\sum_{j \in A_i} I_j - 1\right)I_i,$$

which gives the first term of the bound (4.1). Referring to (3.4), if we take the reference measure $\boldsymbol{\nu}(ds, \{i\}) = \mu_i(ds)$, then $\phi((s,i)) = p_i$ and for $i_1, \ldots, i_k \in \mathcal{I}$, where $i_1, \ldots, i_k$ are all different,

$$j_k((s_1, i_1), \ldots, (s_k, i_k)) = \mathbb{P}(C_{i_1,\ldots,i_k}),$$

where

$$C_{i_1,\ldots,i_k} := \{I_l = 1 \text{ for } l = i_1, \ldots, i_k \text{ and } I_l = 0 \text{ for } l \neq i_1, \ldots, i_k\}.$$

For $\alpha = (s, i)$, $\boldsymbol{\beta} = ((s_1, i_1), \ldots, (s_k, i_k)) \in (A_{(s,i)}^c)^k$, the numerator of (2.9) becomes

$$\sum_{r \geq 0} \frac{1}{r!} \sum_{\{j_1,\ldots,j_r\} \subset A_i \setminus \{i\}} \mathbb{P}(C_{i,i_1,\ldots,i_k,j_1,\ldots,j_r})$$

$$= \mathbb{P}(I_j = 1 \text{ for } j = i, i_1, \ldots, i_k \text{ and } I_j = 0 \text{ for } j \in A_i^c \setminus \{i_1, \ldots, i_k\});$$

and the denominator of (2.9) is reduced to

$$\sum_{r \geq 0} \frac{1}{r!} \sum_{\{j_1,\ldots,j_r\} \subset A_i} \mathbb{P}(C_{i_1,\ldots,i_k,j_1,\ldots,j_r})$$

$$= \mathbb{P}(I_j = 1 \text{ for } j = i_1, \ldots, i_k \text{ and } I_j = 0 \text{ for } j \in A_i^c \setminus \{i_1, \ldots, i_k\}).$$

It follows that

$$\mathcal{G}((s, i), ((s_1, i_1), \ldots, (s_k, i_k)))$$

$$= \mathbb{P}(I_i = 1 | I_j = 1 \text{ for } j = i_1, \ldots, i_k \text{ and } I_j = 0 \text{ for } j \in A_i^c \setminus \{i_1, \ldots, i_k\}).$$

Therefore,

$$\mathcal{G}((s, i), \mathcal{M}_l^{((s,i))}) = \mathbb{E}(I_i | I_j; j \notin A_i)$$

and

$$\int_{(s,i) \in \Gamma} \mathbb{E}|\mathcal{G}((s, i), \mathcal{M}_l^{((s,i))}) - \phi((s, i))|\boldsymbol{\nu}(ds, \{i\}) = \sum_{i \in \mathcal{I}} \mathbb{E}|\mathbb{E}(I_i | I_j, j \notin A_i) - p_i|,$$

which gives $\varepsilon_1$ of Theorem 4.1. On the other hand, in view of $\varepsilon_2$ in (3.5), we can write the Palm process associated with $\mathcal{M}_l$ at $(s, i)$ as

$$\mathcal{M}_{(s,i)}(dt, \{j\}) = \begin{cases} J_{ji}\delta_{\mathcal{U}_j}(dt), & \text{if } j \neq i, \\ 0, & \text{if } j = i \text{ and } t \neq s, \\ \delta_t(dt), & \text{if } j = i \text{ and } t = s. \end{cases}$$

With this coupling, we have $|(\mathcal{M}_{(s,i)})^{((s,i))}| = \sum_{j \notin A_i} J_{ji}$ and $|\mathcal{M}_l^{((s,i))}| = V_i$. So,

$$\frac{1}{|(\mathcal{M}_{(s,i)})^{((s,i))}| \wedge |\mathcal{M}_l^{((s,i))}| + 1} = \frac{1}{V_i \wedge \sum_{j \notin A_i} J_{ji} + 1}$$

and

$$\rho_1''((\mathcal{M}_{(s,i)})^{((s,i))}, \mathcal{M}_l^{((s,i))}) \leq \sum_{j \notin A_i} |J_{ji} - I_j|,$$

which yields $\varepsilon_2$ of Theorem 4.1. Finally, since

$$\boldsymbol{\lambda}(ds, \{i\})\boldsymbol{\lambda}(dt, \{j\}) = p_i p_j \mu_i(ds)\mu_j(dt),$$

the last term of (4.1) follows from the last term of (3.3). $\square$

Bounds on $\mathbb{E}[\frac{1}{V_i+1}|I_j = 1]$ and $\mathbb{E}[\frac{1}{V_i+1}|I_j = I_i = 1]$ may be obtained by applying Lemma 3.5. Sharper bounds can be achieved if additional information about the relationship of $I_i$'s is available, for example, if $I_i$'s are independent.

REMARK 4.4. If $I_i$, $i \in \mathcal{I}$, are locally dependent with neighborhoods $(A_i; i \in \mathcal{I})$, then

$$\mathbb{E}\left[\frac{1}{V_i + 1}\Big|I_j = 1\right] \le \mathbb{E}\left[\frac{1}{V_{ij} + 1}\right],$$

where $V_{ij} = \sum_{k \notin A_i \cup A_j} I_k$.

Random indicators $(I_j; i \in \mathcal{I})$ are said to be *negatively related* (resp. *positively related*) if, for each $i$, $(J_{ji}, j \in \mathcal{I})$ can be constructed in such a way that $J_{ji} \le$ (resp. $\ge$) $I_j$ for $j \in \mathcal{I}$, $j \ne i$ [see Barbour, Holst and Janson (1992), page 24].

PROPOSITION 4.5. *Suppose $(I_j; j \in \mathcal{I})$ are negatively related, and let $\lambda = \mathbb{E}\sum_{i \in \mathcal{I}} I_i$; then*

$$\mathbb{E}\frac{1}{\sum_{i \in \mathcal{I}} I_i + 1} \le \frac{1 - e^{-\lambda}}{\lambda}.$$

PROOF. Indeed, since $(I_j; j \in \mathcal{I})$ are negatively related, for decreasing function $\Phi$,

$$\mathbb{E}\left(\Phi\left(\sum_{i \in \mathcal{I}\setminus\{j\}} I_i\right)\Big|I_j = 1\right) \ge \mathbb{E}\left(\Phi\left(\sum_{i \in \mathcal{I}\setminus\{j\}} I_i\right)\Big|I_j = 0\right),$$

so for fixed $0 < z < 1$, $\mathbb{E}(z^{\sum_{i \in \mathcal{I}\setminus\{j\}} I_i}|I_j)$ is increasing in $I_j$ and $z^{I_j}$ is a decreasing function in $I_j$, giving

$$\mathbb{E}z^{\sum_{i \in \mathcal{I}} I_i} = \mathbb{E}[\mathbb{E}(z^{\sum_{i \in \mathcal{I}\setminus\{j\}} I_i}|I_j)z^{I_j}]$$

$$\le \mathbb{E}[\mathbb{E}(z^{\sum_{i \in \mathcal{I}\setminus\{j\}} I_i}|I_j)]\mathbb{E}[z^{I_j}] = \mathbb{E}(z^{\sum_{i \in \mathcal{I}\setminus\{j\}} I_i})\mathbb{E}z^{I_j},$$

[see Liggett (1985), page 78]. Since $\mathcal{I}$ is a finite or infinitely countable index set, by mathematical induction,

$$\mathbb{E}z^{\sum_{i \in \mathcal{I}} I_i} \le \prod_{i \in \mathcal{I}} \mathbb{E}z^{I_i}.$$

Hence

$$\mathbb{E}\frac{1}{\sum_{i\in\mathcal{I}}I_i+1}=\mathbb{E}\int_0^1 z^{\sum_{i\in\mathcal{I}}I_i}\,dz\leq\int_0^1\prod_{i\in\mathcal{I}}(1-p_i+p_iz)\,dz$$

$$\leq\int_0^1\prod_{i\in\mathcal{I}}e^{-p_i(1-z)}\,dz=\frac{1-e^{-\lambda}}{\lambda}.\qquad\square$$

COROLLARY 4.6. *With the same setup as in Theorem* 4.1, *suppose* $(I_j;j\in\mathcal{I})$ *are negatively related; then*

$$(4.2)\quad d_2(\mathcal{LM},\mathrm{Po}(\boldsymbol{\lambda}))\leq\mathbb{E}\sum_{i\in\mathcal{I}}\left(\frac{5}{\lambda}+\frac{3}{\sum_{j\neq i}J_{ji}+1}\right)\left[p_i^2+p_i\sum_{j\neq i}[I_j-J_{ji}]\right].$$

PROOF. By Theorem 4.1 with $A_i=\{i\}$ and $\varepsilon_2$, the first term of (4.1) vanishes and the last two terms of (4.1) can be rewritten as (4.2). $\square$

As we need to bound $\mathbb{E}[(V_i+1)^{-1}|I_i=I_j=1]$, it is relevant to ask whether $(J_{ki},k\in\mathcal{I})$ are also negatively (resp. positively) related if $(I_j;j\in\mathcal{I})$ are negatively (resp. positively) related. The answer is generally negative, as the following counterexample shows.

COUNTEREXAMPLE 4.7. Choose four sets $B_i$, $1\leq i\leq 4$, so that $\mathbb{P}(B_i)=q$, $\mathbb{P}(B_iB_j)=bq^2$, $\mathbb{P}(B_iB_jB_k)=bq^3$, for all different $1\leq i,j,k\leq 4$; and $\mathbb{P}(B_1B_2B_3B_4)=bq^4$ with $b\leq 2$ and $q$ sufficiently small (e.g., $\leq 0.01$) so that the sets are properly defined. Set $I_i=\mathbb{1}_{B_i}$. Then for any increasing function $\Phi$ on $\{0,1\}^3$ [see Barbour, Holst and Janson (1992), page 27], we have

$$\mathbb{E}[\Phi(I_1,I_2,I_3)|I_4=1]-\mathbb{E}\Phi(I_1,I_2,I_3)$$

$$=q(b-1)[\Phi(1,0,0)+\Phi(0,1,0)+\Phi(0,0,1)-3\Phi(0,0,0)].$$

Hence, by Theorem 2.D of Barbour, Holst and Janson (1992), if we choose $b>$ (resp. $<$) 1, then $(I_i;1\leq i\leq 4)$ are positively (resp. negatively) related. But

$$\mathbb{P}(J_{31}=J_{41}=1|J_{21}=1)=\mathbb{P}(I_3=I_4=1|I_1=I_2=1)=q^2$$

and

$$\mathbb{P}(J_{31}=J_{41}=1)=\mathbb{P}(I_3=I_4=1|I_1=1)=bq^2,$$

so

$$\mathbb{P}(J_{31}=J_{41}=1|J_{21}=1)<(\text{resp. }>)\,\mathbb{P}(J_{31}=J_{41}=1),$$

which implies that $(J_{k1},k=1,\dots,4)$ are not positively (resp. negatively) related.

**5. Applications.** In this section, we apply the main results in Sections 3 and 4 to the Matérn hard-core process, an occupancy problem and rare words in DNA sequences, all of which are different in nature. The results in Section 4 can also be applied to random graphs, for example, to the isolated vertices resulting from the deletion with small probability of each of the edges of a connected graph, where the resulting isolated vertices may remain in their original positions or may be distributed independently and randomly in a carrier space. Since this random graph problem is similar in nature to that of rare words in DNA sequences, it will not be discussed further in this section. A special case of this problem which involves counting the number of isolated vertices has been considered by Roos (1994) and Eichelsbacher and Roos (1999).

5.1. *Matérn hard-core process.* Consider a Poisson number, with mean $\mu$, of points placed independently and uniformly at random in $\Gamma$, where $\Gamma$ is a compact subset of $\mathbb{R}^d$ with volume $V(\Gamma) \neq 0$. A Matérn hard-core process $\Xi$ is produced by deleting any point within distance $r$ of another point, irrespective of whether the latter point has itself already been deleted [see Cox and Isham (1980), page 170]. More precisely, let $\{\alpha'_n\}$ be a realization of points of the Poisson process. Then the points deleted are

$$\{\alpha''_n\} = \{x \in \{\alpha'_n\} : |x - y| < r \text{ for some } y \neq x, y \in \{\alpha'_n\}\},$$

and $\{\alpha_n\} := \{\alpha'_n\} \setminus \{\alpha''_n\}$ constitutes a realization of the Matérn hard-core process $\Xi$ [see Daley and Vere-Jones (1988)].

The Matérn hard-core process is one of the hard-core processes introduced in statistical mechanics to model the distribution of particles with repulsive interactions [see Ruelle (1969), page 6]. It is a special case of the distance models [see Matérn (1986), page 37] and is also a model for underdispersion [see Daley and Vere-Jones (1988), page 366].

Let $X_1, X_2, \ldots$ be independent uniform random variables on $\Gamma$, and let $N$ be a Poisson random variable with mean $\mu$ and independent of $\{X_i; i \geq 1\}$. Then the Poisson process for the arrival points in $\Gamma$ is $Z = \sum_{i=1}^N \delta_{X_i}$. Let $B(x, r) = \{y \in \Gamma : 0 < d_0(y, x) < r\}$, the $r$-neighborhood of $x$, where $d_0(x, y) = |x - y| \wedge 1$. Then the Matérn hard-core process $\Xi$ can be written as $\Xi = \sum_{i=1}^N \delta_{X_i} \mathbb{1}_{\{Z(B(X_i, r))=0\}}$. Also,

$$\Xi(d\alpha) = \sum_{i=1}^N \delta_{X_i}(d\alpha) \mathbb{1}_{\{Z(B(X_i, r))=0\}} = \mathbb{1}_{\{Z(B(\alpha, r))=0\}} Z(d\alpha).$$

Let $\kappa_d$ be the volume of the unit ball in $\mathbb{R}^d$ and let $d_2$ be the second Wasserstein metric generated from $d_0$ as in Section 3.

THEOREM 5.1. *The mean measure of $\Xi$ is $\boldsymbol{\lambda}(d\alpha) = e^{-\mu V(\alpha,r)/V(\Gamma)}\mu \times V(\Gamma)^{-1} d\alpha$, and*

$$d_2(\mathcal{L}\Xi, \text{Po}(\boldsymbol{\lambda})) \leq 10\vartheta + 6\vartheta[3 + (1 - e^{-2^{-d}\vartheta})\vartheta]/(1 + (1-2\vartheta)/\lambda),$$

*where $V(\alpha,r)$ is the volume of $B(\alpha,r)$ and $\vartheta = \mu\kappa_d(2r)^d/V(\Gamma)$.*

PROOF. The Poisson property of $Z$ implies that the counts of points in disjoint sets are independent. So

$$\boldsymbol{\lambda}(d\alpha) = \mathbb{E}(\Xi(d\alpha)) = \mathbb{E}\mathbb{1}_{\{Z(B(\alpha,r))=0\}}\mathbb{E}Z(d\alpha) = e^{-\mu V(\alpha,r)/V(\Gamma)}\mu V(\Gamma)^{-1} d\alpha.$$

Also, whether a point outside $B(\alpha,2r) \cup \{\alpha\}$ is deleted or not is independent of the behavior of $Z$ in $B(\alpha,r) \cup \{\alpha\}$. Hence, we choose $A_\alpha = B(\alpha,2r) \cup \{\alpha\}$ so that $\Xi$ is locally dependent with neighborhoods $(A_\alpha; \alpha \in \Gamma)$ and

$$\mathbb{E}\frac{1}{|\Xi^{(\alpha\beta)}| + 1}\Xi(d\alpha)\Xi(d\beta) = \mathbb{E}\frac{1}{|\Xi^{(\alpha\beta)}| + 1}\mathbb{E}\Xi(d\alpha)\Xi(d\beta).$$

Applying Corollary 3.6 gives

$$\begin{aligned}
d_2(\mathcal{L}\Xi, \text{Po}(\boldsymbol{\lambda})) &\leq \int_{\alpha \in \Gamma}\int_{\beta \in A_\alpha \setminus \{\alpha\}} \left(\frac{5}{\lambda} + \mathbb{E}\frac{3}{|\Xi^{(\alpha\beta)}| + 1}\right)\mathbb{E}\Xi(d\alpha)\Xi(d\beta) \\
&\quad + \int_{\alpha \in \Gamma}\int_{\beta \in A_\alpha} \left(\frac{5}{\lambda} + \mathbb{E}\frac{3}{|\Xi^{(\alpha\beta)}| + 1}\right)\boldsymbol{\lambda}(d\alpha)\boldsymbol{\lambda}(d\beta).
\end{aligned}$$ (5.1)

Now,

$$\mathbb{E}|\Xi^{(\alpha\beta)}| = \int_{\Gamma_{\alpha\beta}} e^{-\mu_\Gamma V(x,r)}\mu_\Gamma \, dx,$$

where $\Gamma_{\alpha\beta} = \Gamma \setminus (A_\alpha \cup A_\beta)$ and $\mu_\Gamma = \mu/V(\Gamma)$. On the other hand,

$\mathbb{E}\Xi(d\alpha)\Xi(d\beta)$

$$= \begin{cases}
e^{-\mu_\Gamma(V(\alpha,r)+V(\beta,r))}\mu_\Gamma^2 \, d\alpha \, d\beta, & \text{if } |\alpha - \beta| \geq 2r, \\
e^{-\mu_\Gamma(V(\alpha,r)+V(\beta,r)-V(\alpha,\beta,r))}\mu_\Gamma^2 \, d\alpha \, d\beta, & \text{if } r \leq |\alpha - \beta| < 2r, \\
0, & \text{if } 0 < |\alpha - \beta| < r, \\
e^{-\mu_\Gamma V(\alpha,r)}\mu_\Gamma \, d\alpha, & \text{if } \alpha = \beta,
\end{cases}$$

where $V(\alpha,\beta,r)$ is the volume of $B(\alpha,r) \cap B(\beta,r)$. Hence,

$$\begin{aligned}
\mathbb{E}[|\Xi^{(\alpha\beta)}|^2] &= \mathbb{E}\iint_{x,y \in \Gamma_{\alpha\beta}} \Xi(dx)\Xi(dy) \\
&= \int_{\Gamma_{\alpha\beta}} e^{-\mu_\Gamma V(x,r)}\mu_\Gamma \, dx \\
&\quad + \iint_{x,y \in \Gamma_{\alpha\beta}, |x-y| \geq 2r} e^{-\mu_\Gamma(V(x,r)+V(y,r))}\mu_\Gamma^2 \, dx \, dy \\
&\quad + \iint_{x,y \in \Gamma_{\alpha\beta}, r \leq |x-y| < 2r} e^{-\mu_\Gamma(V(x,r)+V(y,r)-V(x,y,r))}\mu_\Gamma^2 \, dx \, dy.
\end{aligned}$$

Writing

$$[\mathbb{E}|\Xi^{(\alpha\beta)}|]^2 = \iint_{x,y \in \Gamma_{\alpha\beta}} e^{-\mu_\Gamma(V(x,r)+V(y,r))} \mu_\Gamma^2 \, dx \, dy,$$

we have

$$\begin{aligned}
\operatorname{Var}(|\Xi^{(\alpha\beta)}|) &= \int_{\Gamma_{\alpha\beta}} e^{-\mu_\Gamma V(x,r)} \mu_\Gamma \, dx \\
&\quad + \iint_{x,y \in \Gamma_{\alpha\beta}, r \le |x-y| < 2r} e^{-\mu_\Gamma(V(x,r)+V(y,r)-V(x,y,r))} \mu_\Gamma^2 \, dx \, dy \\
&\quad - \iint_{x,y \in \Gamma_{\alpha\beta}, |x-y| < 2r} e^{-\mu_\Gamma(V(x,r)+V(y,r))} \mu_\Gamma^2 \, dx \, dy \\
&\le \int_{\Gamma_{\alpha\beta}} e^{-\mu_\Gamma V(x,r)} \mu_\Gamma \, dx \\
&\quad + \iint_{x,y \in \Gamma_{\alpha\beta}, |x-y| < 2r} e^{-\mu_\Gamma V(x,r)} [1 - e^{-\mu_\Gamma V(y,r)}] \mu_\Gamma^2 \, dx \, dy \\
&\le \{1 + (1 - e^{-\mu_\Gamma \kappa_d r^d}) \mu_\Gamma \kappa_d (2r)^d\} \int_{\Gamma_{\alpha\beta}} e^{-\mu_\Gamma V(x,r)} \mu_\Gamma \, dx.
\end{aligned}$$

Thus,

$$\kappa = \frac{\operatorname{Var}(|\Xi^{(\alpha\beta)}|+1)}{\mathbb{E}(|\Xi^{(\alpha\beta)}|+1)} \le \frac{\operatorname{Var}(|\Xi^{(\alpha\beta)}|)}{\mathbb{E}(|\Xi^{(\alpha\beta)}|)} \le 1 + (1 - e^{-\mu_\Gamma \kappa_d r^d}) \mu_\Gamma \kappa_d (2r)^d,$$

which, together with Lemma 3.5, yields

$$\mathbb{E} \frac{1}{|\Xi^{(\alpha\beta)}|+1} \le \frac{2+\kappa}{\int_{\Gamma_{\alpha\beta}} e^{-\mu_\Gamma V(x,r)} \mu_\Gamma \, dx + 1} \le \frac{3 + (1 - e^{-\mu_\Gamma \kappa_d r^d}) \mu_\Gamma \kappa_d (2r)^d}{\lambda + 1 - 2\mu_\Gamma \kappa_d (2r)^d}.$$

Finally,

$$\int_{\alpha \in \Gamma} \int_{\beta \in A_\alpha \setminus \{\alpha\}} \mathbb{E}\Xi(d\alpha)\Xi(d\beta) \le \int_{\alpha \in \Gamma} \int_{\beta \in A_\alpha} e^{-\mu_\Gamma V(\alpha,r)} \mu_\Gamma^2 \, d\alpha \, d\beta \le \mu_\Gamma \kappa_d (2r)^d \lambda$$

and

$$\int_{\alpha \in \Gamma} \int_{\beta \in A_\alpha} \boldsymbol{\lambda}(d\alpha)\boldsymbol{\lambda}(d\beta) \le \mu_\Gamma \kappa_d (2r)^d \lambda.$$

Applying these inequalities to the relevant terms in (5.1) gives Theorem 5.1. $\square$

5.2. *Occupancy problem.* Suppose $s$ balls are dropped independently into $n$ urns with probability $p_k$ of going into the $k$th urn. Two cases of the distribution of urns with given content have been studied in the literature.

They are urns with at most $m$ balls (*right-hand domain*) and urns with at least $m$ balls (*left-hand domain*), where $m$ is a fixed nonnegative integer [see Kolchin, Sevast'yanov and Chistyakov (1978) and also Barbour, Holst and Janson (1992), Chapter 6]. In this section, we consider the right-hand domain. So far, the focus in the literature has been on the total number of urns [see Arratia, Goldstein and Gordon (1989) and Barbour, Holst and Janson (1992), and references therein] and little attention has been paid to the locations of the urns.

We assume the urns are numbered from 1 to $n$ and let $X_i$ be the number of balls in the $i$th box, $1 \leq i \leq n$. Define a point process $\Xi$ on $\Gamma = [0, 1]$ as follows:

$$\Xi = \sum_{i=1}^{n} \mathbb{1}_{\{X_i \leq m\}} \delta_{i/n}.$$

The mean measure of $\Xi$ is then $\boldsymbol{\mu} = \sum_{i=1}^{n} \pi_i \delta_{i/n}$, where $\pi_i = \sum_{j=0}^{m} \binom{s}{j} p_i^j \times (1 - p_i)^{s-j}$ and $\mu = \sum_{i=1}^{n} \pi_i$. Set $\boldsymbol{\lambda}(dt) = n\pi_i \, dt$ for $(i-1)/n < t \leq i/n$, $i = 1, 2, \ldots, n$ and $d_0(t_1, t_2) = |t_1 - t_2|$ for $t_1, t_2 \in \Gamma$. Let

$$\mu' = \min_{i \neq j, 1 \leq i, j \leq n} \sum_{k \neq i, j, 1 \leq k \leq n} \mathbb{P}(X_k \leq m | X_i = X_j = 0)$$

and

$$\mu'' = \min_{1 \leq i \leq n} \sum_{j \neq i, 1 \leq j \leq n} \mathbb{P}(X_j \leq m | X_i = 0) \geq \mu'.$$

If $s \min_{1 \leq i \leq n} p_i$ is large, then we would expect good Poisson process approximation.

Theorem 5.2. *With the above setup,*

$$d_2(\mathcal{L}\Xi, \mathrm{Po}(\boldsymbol{\lambda}))$$

(5.2)
$$\leq \frac{1}{2n} + \left(\frac{5}{\mu} + \frac{3}{\mu'}\right)[\mathbb{E}(|\Xi|) - \mathrm{Var}(|\Xi|)]$$

(5.3)
$$\leq \frac{1}{2n} + C\left\{\pi_* + \frac{s}{\mu}\left(\frac{\ln s + m \ln \ln s + 5m}{s - \ln s - m \ln \ln s - 4m}\mu + \frac{4}{s}\right)^2\right\},$$

*where* (5.3) *is valid for* $s > \ln s + m \ln \ln s + 4m$,

$$C = 5 + 3\left(\frac{1 - 3p_* + 2p_*^2}{1 - 3p_*}\right)^s \left(1 - \frac{2\pi_*}{\mu}\right)^{-1},$$

*with* $\pi_* = \max_{1 \leq i \leq n} \pi_i$ *and* $p_* = \max_{1 \leq i \leq n} p_i < 1/3$.

PROOF. By the triangle inequality, $d_2(\mathcal{L}\Xi, \mathrm{Po}(\boldsymbol{\lambda})) \leq d_2(\mathcal{L}\Xi, \mathrm{Po}(\boldsymbol{\mu})) + d_2(\mathrm{Po}(\boldsymbol{\mu}), \mathrm{Po}(\boldsymbol{\lambda}))$, so the term $1/(2n)$ follows immediately from estimating $d_2(\mathrm{Po}(\boldsymbol{\mu}), \mathrm{Po}(\boldsymbol{\lambda}))$ [see Brown and Xia (1995), (2.8)]. For each $1 \leq i \leq n$, let $I_i = \mathbb{1}_{\{X_i \leq m\}}$, then $(I_i; 1 \leq i \leq n)$ are negatively related. Indeed, if $X_i \leq m$, take $Y_{ji} = X_j$ for all $j$. If $X_i > m$, take a random variable $\tilde{X}_i$ which is independent of $\{X_1, \ldots, X_n\}$ and has distribution $\mathcal{L}(X_i | X_i \leq m)$ and take $X_i - \tilde{X}_i$ balls from urn $i$ and redistribute them to the other urns with probabilities $p_j/(1 - p_i)$ for $j \neq i$. Let $Y_{ji}$ be the number of balls in urn $j$ after the redistribution and set $J_{ji} = \mathbb{1}_{\{Y_{ji} \leq m\}}$. This coupling $(J_{ji}; 1 \leq j \leq n)$ satisfies

$$\mathcal{L}(J_{ji}; 1 \leq j \leq n) = \mathcal{L}(I_j; 1 \leq j \leq n | I_i = 1), \qquad J_{ji} \leq I_j \text{ for all } j \neq i$$

[see Barbour, Holst and Janson (1992), page 122].

We have from Corollary 4.6 that

$$(5.4) \quad d_2(\mathcal{L}\Xi, \mathrm{Po}(\boldsymbol{\mu})) \leq \mathbb{E} \sum_{i=1}^{n} \left( \frac{5}{\mu} + \frac{3}{\sum_{j \neq i} J_{ji} + 1} \right) \left( \sum_{k \neq i} (I_k - J_{ki}) \pi_i + \pi_i^2 \right).$$

Now, the above coupling can be modified to show that, for $l \geq 1$,

$$(I_j; j \neq i_1, \ldots, i_l | X_{i_1} = \cdots = X_{i_l} = 0)$$

are also negatively related. In fact, denote $\mathbf{i} = (i_1, \ldots, i_l)$. If $X_{i_1} = \cdots = X_{i_l} = 0$, take $Z'_{j\mathbf{i}} = X_j$ for all $j \neq i_1, \ldots, i_l$. If one of $X_{i_1}, \ldots, X_{i_l}$ is not 0, take all balls in urns $i_1, \ldots, i_l$ and relocate them to the other urns with probabilities $p'_j := p_j/(1 - p_{i_1} - \cdots - p_{i_l})$ for $j \neq i_1, \ldots, i_l$. After the relocation, let $Z'_{j\mathbf{i}}$ be the number of balls in urn $j$ and $J'_{j\mathbf{i}} = \mathbb{1}_{\{Z'_{j\mathbf{i}} \leq m\}}$. Next, for $k \neq i_1, \ldots, i_l$, if $Z'_{k\mathbf{i}} \leq m$, take $J''_{jk\mathbf{i}} = J'_{j\mathbf{i}}$. If $Z'_{k\mathbf{i}} > m$, take a random variable $\tilde{Z}'_{k\mathbf{i}}$ which is independent of $\{Z'_{1\mathbf{i}}, \ldots, Z'_{n\mathbf{i}}\}$ and has distribution $\mathcal{L}(Z'_{k\mathbf{i}} | Z'_{k\mathbf{i}} \leq m)$. Remove $Z'_{k\mathbf{i}} - \tilde{Z}'_{k\mathbf{i}}$ balls from urn $k$ and redistribute them to the other urns with probabilities $p'_j/(1 - p'_k)$, $j \neq k, i_1, \ldots, i_l$. After this, let $J''_{jk\mathbf{i}} = 1$ if there are at most $m$ balls in urn $j$; otherwise, let $J''_{jk\mathbf{i}} = 0$. These couplings satisfy

$$\mathcal{L}(J'_{j\mathbf{i}}; j \neq i_1, \ldots, i_l) = \mathcal{L}(I_j; j \neq i_1, \ldots, i_l | X_{i_1} = \cdots = X_{i_l} = 0),$$

$$\mathcal{L}(J''_{jk\mathbf{i}}; j \neq i_1, \ldots, i_l) = \mathcal{L}(J'_{j\mathbf{i}}; j \neq i_1, \ldots, i_l | J'_{k\mathbf{i}} = 1),$$

$$J''_{jk\mathbf{i}} \leq J'_{j\mathbf{i}} \qquad \text{for all } j \neq k, i_1, \ldots, i_l,$$

$$J'_{j\mathbf{i}} \leq I_j \qquad \text{for all } j \neq i_1, \ldots, i_l.$$

In particular, if $\mathbf{i} = i$, then $J'_{ji} \leq J_{ji}$ for $j \neq i$.

By these couplings and Proposition 4.5, we have

$$(5.5) \qquad \mathbb{E} \frac{1}{\sum_{j \neq i} J_{ji} + 1} \leq \mathbb{E} \frac{1}{\sum_{j \neq i} J'_{ji} + 1} \leq \frac{1}{\sum_{j \neq i} \mathbb{E} J'_{ji}} \leq \frac{1}{\mu'}.$$

On the other hand, for $k \neq i$, denote $\mathbf{k} = (i, k)$, we have

$$\mathbb{E} \frac{I_k - J_{ki}}{\sum_{j \neq i} J_{ji} + 1}$$

$$\leq \mathbb{E} \frac{I_k - J_{ki}}{\sum_{j \neq i, k} J_{ji} + 1}$$

$$= \sum_{l_1=0}^{m} \sum_{l_2=m+1}^{s} \mathbb{E} \left[ \frac{1}{\sum_{j \neq i, k} \mathbb{1}_{\{Y_{ji} \leq m\}} + 1} \Big| X_k = l_1, Y_{ki} = l_2 \right] \mathbb{P}(X_k = l_1, Y_{ki} = l_2)$$

$$\leq \mathbb{E} \frac{1}{\sum_{j \neq i, k} \mathbb{1}_{\{Z'_{j\mathbf{k}} \leq m\}} + 1} \sum_{l_1=0}^{m} \sum_{l_2=m+1}^{s} \mathbb{P}(X_k = l_1, Y_{ki} = l_2)$$

$$= \mathbb{E} \frac{1}{\sum_{j \neq i, k} J'_{j\mathbf{k}} + 1} \mathbb{P}(I_k = 1, J_{ki} = 0)$$

$$\leq \left[ \sum_{j \neq i, k} \mathbb{E} J'_{j\mathbf{k}} \right]^{-1} \mathbb{P}(I_k = 1, J_{ki} = 0)$$

$$\leq \frac{\mathbb{P}(I_k = 1, J_{ki} = 0)}{\mu'} = \frac{\mathbb{E}(I_k - J_{ki})}{\mu'}.$$

Hence,

$$\mathbb{E} \left( \frac{\sum_{k \neq i} I_k - \sum_{k \neq i} J_{ki}}{\sum_{j \neq i} J_{ji} + 1} \right) \leq \frac{1}{\mu'} \mathbb{E} \left[ \sum_{k \neq i} (I_k - J_{ki}) \right],$$

which, combined with (5.4) and (5.5), yields

$$d_2(\mathcal{L}\Xi, \mathrm{Po}(\boldsymbol{\mu})) \leq \left( \frac{5}{\mu} + \frac{3}{\mu'} \right) \mathbb{E} \sum_{i=1}^{n} \left[ \left( \sum_{k \neq i} I_k - \sum_{k \neq i} J_{ki} \right) \pi_i + \pi_i^2 \right].$$

On the other hand, since for $k \neq i$, $\mathbb{E}(J_{ki}) \pi_i = \mathbb{P}(I_k = I_i = 1) = \mathbb{E}(I_k I_i)$, we have

$$\mathbb{E} \sum_{i=1}^{n} \left[ \left( \sum_{k \neq i} I_k - \sum_{k \neq i} J_{ki} \right) \pi_i + \pi_i^2 \right]$$

$$= \sum_{1 \leq i, k \leq n} \mathbb{E}(I_k) \mathbb{E}(I_i) - \sum_{i \neq k, 1 \leq i, k \leq n} \mathbb{E}(I_i I_k)$$

$$= \mathbb{E}(|\Xi|) - \mathrm{Var}(|\Xi|).$$

Therefore, (5.2) follows. To prove (5.3), we note from Theorem 6.D of Barbour, Holst and Jonson [(1992), page 122] that

$$1 - \frac{\mathrm{Var}(|\Xi|)}{\mathbb{E}(|\Xi|)} \leq \pi_* + \frac{s}{\mu} \left( \frac{\ln s + m \ln \ln s + 5m}{s - \ln s - m \ln \ln s - 4m} \mu + \frac{4}{s} \right)^2.$$

So, it remains to show that

$$(5.6) \qquad \frac{\mu}{\mu'} \leq \left(\frac{1 - 3p_* + 2p_*^2}{1 - 3p_*}\right)^s \left(1 - \frac{2\pi_*}{\mu}\right)^{-1}.$$

To prove (5.6), notice that, for $1 \leq i, j \leq n$ with $i \neq j$,

$$\sum_{k \neq i, j} \mathbb{P}(X_k \leq m | X_i = X_j = 0)$$

$$= \sum_{k \neq i, j} \sum_{0 \leq l \leq m} \binom{s}{l} \left(\frac{p_k}{1 - p_i - p_j}\right)^l \left(1 - \frac{p_k}{1 - p_i - p_j}\right)^{s-l}$$

$$\geq \sum_{k \neq i, j} \sum_{0 \leq l \leq m} \binom{s}{l} p_k^l (1 - p_k)^{s-l} \left(\frac{1 - p_i - p_j - p_k}{(1 - p_k)(1 - p_i - p_j)}\right)^s$$

$$\geq \left(\frac{1 - 3p_*}{1 - 3p_* + 2p_*^2}\right)^s \sum_{k \neq i, j} \sum_{0 \leq l \leq m} \binom{s}{l} p_k^l (1 - p_k)^{s-l}$$

$$\geq \left(\frac{1 - 3p_*}{1 - 3p_* + 2p_*^2}\right)^s (\mu - 2\pi_*).$$

Hence

$$\mu' \geq \left(\frac{1 - 3p_*}{1 - 3p_* + 2p_*^2}\right)^s (\mu - 2\pi_*),$$

which implies (5.6).  □

5.3. *Rare words in biomolecular sequences.* One of the important problems in biomolecular sequence analysis is the study of the distribution of words in a DNA sequence. A DNA sequence may be regarded as a sequence of letters taken from the alphabet {A, C, G, T}. The letters A, C, G, T represent the four nucleotides: adenine, cytosine, guanine and thymine, respectively. They form two complementary pairs, namely {A, T} and {C, G}.

It is known that repetition of a given word or a group of words or occurrences of unusually large clusters of words are known to have biological functions. For example, unusually large clusters of palindromes are known to contain such significant sites as origins of replication and gene regulators. Here palindromes are symmetrical words of DNA in the sense that they read exactly the same as their reverse complementary sequences. In Leung and Yamashita (1999), palindromes of certain lengths are assumed to be independent and uniformly distributed in herpesvirus genomes, and the *r*-scan statistic is used to identify unusually large clusters of palindromes.

It is commonly assumed that the bases of DNA are independent random variables taking values in the set {A, C, G, T}. Under this assumption, if each word of a particular type is represented by a point, then the points representing these words form a locally dependent point process. Theorem 4.1 in this paper provides an error bound for approximating such a point process by a Poisson process. The error bound can be used to find conditions for which the approximation is good. In general, the approximation is good if the words are rare in the sense that the probabilities of their occurrences are small. However, the error bound can be made more explicit only when the words are specified.

As an application, Theorem 4.1 provides a mathematical basis for Poisson process modeling of rare words in a biomolecular sequence, and in particular of palindromes in a DNA sequence. A consequence of this is that the observed rare words may be regarded as a realization of i.i.d. random variables, thus providing a mathematical basis for the assumption in Leung and Yamashita (1999) that the points representing the palindromes are independent and uniformly distributed in the herpesvirus genomes.

In Leung, Xia and Chen (2002) Poisson process approximation for palindromes in sixteen herpesvirus genomes is studied. The centers of palindromes in each herpesvirus genome are represented by the point process on $\{0, 1/n, 2/n, \ldots, (n-1)/n, 1\}$:

$$\Xi = \sum_{i=1}^{n} I_i \delta_{i/n}, \tag{5.7}$$

where the length of genome (number of base pairs) is denoted by $M$, those palindromes considered are of length at least $2L$ (the length must be even) and called $2L$-palindromes, the center of a palindrome of length $2K$ is the $K$th base in the palindrome from the left, and the number of possible centers of $2L$-palindromes is $M - 2L + 1$, denoted by $n$. Also, $I_i$ is the indicator random variable for the occurrence of a $2L$-palindrome centered at base $i + L - 1$ of the DNA sequence. The palindromes are represented by their centers because the latter are fixed irrespective of the lengths of the former, whereas the first base pair of a palindrome of at least a certain length is random and will give rise to complications in the analysis if it is used to represent the palindrome.

Since $2L$-palindromes with centers sufficiently far apart (more specifically, further than $2L - 1$ bases apart) occur independently, the point process (5.7) is a sequence of marked locally dependent trials as described in Section 4 of this paper, to which Theorem 4.1 is applicable. Here $(I_i; 1 \leq i \leq n)$ are locally dependent with neighborhoods

$$A_i = \{j : i - 2L + 1 \leq j \leq i + 2L - 1\} \cap \{1, 2, \ldots, n\}, \qquad i = 1, 2, \ldots, n.$$

Take $\Gamma = [0, 1]$ and $d_0(x, y) = |x - y|$. Let $p_i = \mathbb{P}(I_i = 1)$ and $p_{ij} = \mathbb{P}(I_i = I_j = 1)$. It can be shown that $p_i = \theta^L$, where $\theta = 2(p_A p_T + p_C p_G)$ and $p_A$, $p_T$, $p_C$, $p_G$ are the probabilities of A, T, C, G, respectively.

Suppose

$$p_A = p_T, \qquad p_C = p_G \quad \text{and} \quad 4 \leq L \leq \frac{n}{500}. \tag{5.8}$$

Then the next theorem follows from Theorem 4.1 with $\mathcal{U}_i = i/n$, Lemma 3.5 and a two-step approximation as in Section 5.2; namely, first approximate $\Xi$ by a Poisson process with the same mean measure as that of $\Xi$ and then approximate the latter by a Poisson process on $[0, 1]$ with intensity measure $\lambda \, dx$.

THEOREM 5.3. *We have*

$$d_2(\mathcal{L}\Xi, \mathrm{Po}(\boldsymbol{\lambda})) \leq \frac{26}{\lambda}(b_1 + b_2) + \frac{1}{2n} \leq 131 L \theta^{L/2}, \tag{5.9}$$

*where* $\lambda = \sum_{i=1}^{n} p_i = n\theta^L$, $b_1 = \sum_{i=1}^{n} \sum_{j \in A_i} p_i p_j \leq n(4L-1)\theta^{2L}$,

$$b_2 = \sum_{i=1}^{n} \sum_{j \in A_i, j \neq i} p_{ij} \leq n(4L-2)\theta^{3L/2}$$

*and* $\boldsymbol{\lambda}(dx) = \lambda \, dx$.

Since a proof of Theorem 5.3 is given in Leung, Choi, Xia and Chen (2002), we will not give one here. It suffices to mention that the explicit bound on the overlap probabilities in (5.9) is due to the explicit nature of the palindrome. In order for $\mathrm{Po}(\boldsymbol{\lambda})$ to be nondegenerate in the limit, $\lambda = n\theta^L$ must converge to a positive number as $n \to \infty$. This means that $L = \ln n / \ln(1/\theta) + d$, where $d$ is bounded. For such an $L$, the assumption (5.8) is satisfied for sufficiently large $n$, Theorem 5.3 holds and the upper bound in (5.9) tends to 0 as $n \to \infty$.

A significant feature of the bound in (5.9) is that it has the Stein factor $1/\lambda$. This is crucial for accuracy, as the value of $\lambda$ ranges from about 100 to 300 for the sixteen herpesvirus genomes under study.

In Leung, Choi, Xia and Chen (2002), a direct proof of a special case of Theorem 4.1 with $\mathcal{U}_i = i/n$ is given (see Theorem 1 and the Appendix). Also given are the details of deducing Theorem 5.3 from the special case of Theorem 4.1 and the proof of the upper bound $131 L \theta^{L/2}$ (see Lemmas 1 and 2 and Propositions 1 and 2). This upper bound is then used as a guide to choose optimal lengths of palindromes for the approximation. The scan statistics is then applied to identify unusually large clusters of palindromes for each of the sixteen herpesviruses.

**Acknowledgment.** The authors thank Peter Hall for this invitation to the Australian National University in June 2003, during which part of the paper was written.

# REFERENCES

ARRATIA, R., GOLDSTEIN, L. and GORDON, L. (1989). Two moments suffice for Poisson approximations: The Chen–Stein method. *Ann. Probab.* **17** 9–25. MR972770

ARRATIA, R., GOLDSTEIN, L. and GORDON, L. (1990). Poisson approximation and the Chen–Stein method. *Statist. Sci.* **5** 403–434. MR1092983

BARBOUR, A. D. (1988). Stein's method and Poisson process convergence. *J. Appl. Probab.* **25A** 175–184. MR974580

BARBOUR, A. D. and BROWN, T. C. (1992). Stein's method and point process approximation. *Stochastic Process. Appl.* **43** 9–31. MR1190904

BARBOUR, A. D., HOLST, L. and JANSON, S. (1992). *Poisson Approximation.* Oxford Univ. Press. MR1163825

BARBOUR, A. D. and MÅNSSON, M. (2002). Compound Poisson process approximation. *Ann. Probab.* **30** 1492–1537. MR1920275

BILLINGSLEY, P. (1968). *Convergence of Probability Measures.* Wiley, New York. MR233396

BROWN, T. C., WEINBERG, G. V. and XIA, A. (2000). Removing logarithms from Poisson process error bounds. *Stochastic Process. Appl.* **87** 149–165. MR1751169

BROWN, T. C. and XIA, A. (1995). On metrics in point process approximation. *Stochastics Stochastics Rep.* **52** 247–263. MR1381671

BROWN, T. C. and XIA, A. (2001). Stein's method and birth–death processes. *Ann. Probab.* **29** 1373–1403. MR1872746

CHEN, L. H. Y. (1975). Poisson approximation for dependent trials. *Ann. Probab.* **3** 534–545. MR428387

CHEN, L. H. Y. (1998). Stein's method: Some perspectives with applications. *Probability Towards 2000. Lecture Notes in Statist.* **128** 97–122. Springer, Berlin. MR1632651

COX, D. R. and ISHAM, V. (1980). *Point Processes.* Chapman and Hall, London. MR598033

DALEY, D. J. and VERE-JONES, D. (1988). *An Introduction to the Theory of Point Processes.* Springer, Berlin. MR950166

EICHELSBACHER, P. and ROOS, M. (1999). Compound Poisson approximation for dissociated random variables via Stein's method. *Combin. Probab. Comput.* **8** 335–346. MR1723647

JANOSSY, L. (1950). On the absorption of a nucleon cascade. *Proc. Roy. Irish Acad. Sci. Sect. A* **53** 181–188. MR45341

KALLENBERG, O. (1983). *Random Measures.* Academic Press, London. MR818219

KOLCHIN, V. F., SEVAST'YANOV, B. A. and CHISTYAKOV, V. P. (1978). *Random Allocations.* Winston, Washington, DC. MR471016

LEUNG, M. Y., CHOI, K. P., XIA, A. and CHEN, L. H. Y. (2002). Nonrandom clusters of palindromes in herpesvirus genomes. Preprint.

LEUNG, M. Y. and YAMASHITA, T. E. (1999). Applications of the scan statistic in DNA sequence analysis. In *Scan Statistics and Applications* 269–286. Birkhäuser, Boston. MR1697750

LIGGETT, T. M. (1985). *Interacting Particle Systems.* Springer, New York. MR776231

MATÉRN, B. (1986). *Spatial Variation,* 2nd ed. *Lecture Notes in Statist.* **36**. Springer, New York. MR867886

NIELSEN, O. A. (1997). *An Introduction to Integration and Measure Theory.* Wiley, New York. MR1468232

PARTHASARATHY, K. R. (1967). *Probability Measures on Metric Spaces.* Academic Press, New York. MR226684

Roos, M. (1994). Stein's method for compound Poisson approximation: The local approach. *Ann. Appl. Probab.* **4** 1177–1187. MR1304780

Ruelle, D. (1969). *Statistical Mechanics*: *Rigorous Results*. Benjamin, New York. MR289084

Srinivasan, S. K. (1969). *Stochastic Theory and Cascade Processes*. North-Holland, Amsterdam. MR261711

Stein, C. (1972). A bound for the error in the normal approximation to the distribution of a sum of dependent random variables. *Proc. Sixth Berkeley Symp. Math. Statist. Probab.* **3** 583–602. Univ. California Press, Berkeley. MR402873

Xia, A. (1997). On using the first difference in the Stein–Chen method. *Ann. Appl. Probab.* **7** 899–916. MR1484790

Xia, A. (2000). Poisson approximation, compensators and coupling. *Stochastic Anal. Appl.* **18** 159–177. MR1739289

Institute for Mathematical Sciences
National University of Singapore
3 Prince George's Park
Singapore 118402
Republic of Singapore
e-mail: lhychen@ims.nus.edu.sg

Department of Mathematics
  and Statistics
University of Melbourne
Victoria 3010
Australia
e-mail: xia@ms.unimelb.edu.au